\documentclass[11pt]{article}

\usepackage{graphicx}

\def\be{\begin{equation}}
\def\ee{\end{equation}}
\def\x{{\bf x} }

\begin{document}

\title{\Large \bf Review of Metaheuristics  and
Generalized Evolutionary Walk Algorithm}

\author{Xin-She Yang \\
Mathematics and Scientific Computing \\ National Physical Laboratory, Teddington, TW11 0LW, UK}

\date{}
\maketitle

\begin{abstract}
Metaheuristic algorithms are often nature-inspired, and they are becoming very
powerful in solving global optimization problems. More than a dozen
of major metaheuristic algorithms have been developed over the last three decades,
and there exist even more variants and hybrid of metaheuristics. This paper intends
to provide an overview of nature-inspired metaheuristic algorithms, from a brief
history to their applications. We try to analyze the main components of these
algorithms and how and why they works. Then, we intend to provide a unified view
of metaheuristics by proposing a generalized evolutionary walk algorithm (GEWA).
Finally, we discuss some of the important open questions.

{\bf Keywords:} Algorithms; ant colony optimization; cuckoo search; differential evolution; firefly algorithm; harmony search;
genetic algorithm; metaheuristic; simulated annealing; particle swarm optimisation. \\

{\bf Reference} to this paper should be made as follows: \\
Xin-She Yang (2011) `Review of metaheuristics and generalized evolutionary walk algorithm',
{\it Int. J. Bio-Inspired Computation}, Vol. 3, No. 2, pp. 77-84 (2011).

\end{abstract}
\newpage

\section{Introduction}

Nature-inspired metaheuristic algorithms are becoming increasingly popular
in optimisation and applications over the last three decades. There are many reasons
for this popularity and success, and one of the main reasons is that these
algorithms have been developed by mimicking the most successful processes in nature,
including biological systems, and physical and chemical processes. For most
algorithms, we know their fundamental components, how exactly they interact to
achieve efficiency still remains partly a mystery, which inspires more active studies.
Convergence analysis of a few algorithms such as the particle swarm optimisation
shows some insight, but in general mathematical analysis of metaheuristic algorithms
remains unsolved and still an ongoing active research topic.

The solution of an optimisation problem requires the choice and the correct use of
the right algorithm. The choice of an algorithm largely depends on the
type of the optimisation problem at hand. For large-scale nonlinear global
optimisation problems, there is no agreed guideline for how to choose and what to choose.
In fact, we are not even sure whether an efficient algorithm exists,
which is especially true for NP-hard problems, and most real-world problem
often are NP-hard indeed.

In most applications, we can in general write an optimisation problem as the following
generic form
\be \textrm{minimise }_{\x \in \Re^n} \;\;\;\; f_i(\x), \qquad (i=1,2,...,M),  \label{Engopt-equ-100} \ee
\be \textrm{subject to } \; h_j(\x) =0, \;\; (j=1,2,...,J),  \ee
\be \qquad \qquad \quad g_k(\x) \le 0, \;\; (k=1,2,...,K), \ee
where $f_i(\x), h_j(\x)$ and $g_k(\x)$ are functions of the design vector
\be \x=(x_1, x_2, ..., x_n)^T. \ee Here the components $x_i$ of $\x$ are called design
or decision variables, and they can be real continuous, discrete or the mixed of these two.
The functions $f_i(\x)$ where $i=1,2,...,M$
are called the objective functions or simply cost functions, and in the case of
$M=1$, there is only a single objective.
The space spanned by the decision variables is called the design space
or  search space $\Re^n$,
while the space formed by the objective function values is called the solution space or
response space.
The equalities for $h_j$ and
inequalities for $g_k$ are called constraints. It is worth pointing out that
we can also write the inequalities in the other way $\ge 0$, and we can also
formulate the objectives as a maximisation problem.

The algorithms used for solving optimisation problems can be very
diverse, from conventional algorithms to modern metaheuristics.

Most conventional or classic algorithms are deterministic.
For example, the simplex method in linear programming is deterministic.
Some deterministic optimisation algorithms used the gradient information,
they are called gradient-based algorithms. For example, the well-known
Newton-Raphson algorithm is gradient-based, as it uses the function
values and their derivatives, and it works extremely well
for smooth unimodal problems. However, if there is some discontinuity
in the objective function, it does not work well. In this case,
a non-gradient algorithm is preferred. Non-gradient-based, or gradient-free/derivative-free,
algorithms do not use any derivative, but only the function values.
Hooke-Jeeves pattern search and Nelder-Mead downhill simplex are
examples of gradient-free algorithms.

For stochastic algorithms, in general we have two types: heuristic
and metaheuristic, though their difference is small. Loosely speaking,
{\it heuristic} means `to find' or `to discover by trial and error'.
Quality solutions to a tough optimisation problem can be found in
a reasonable amount of time, but there is no guarantee that optimal
solutions are reached. It hopes that these algorithms work most
of the time, but not necessarily  all the time. This is good when we do not necessarily want
the best solutions but rather good solutions which are easily reachable.

\section{Metaheuristics}

In  metaheuristic algorithms,  {\it meta-} means  `beyond' or `higher level',
and they generally perform better than simple heuristics.
All metaheuristic algorithms use certain tradeoff of local search and global exploration.
Variety of solutions are often realized via randomisation.
Despite the popularity of metaheuristics, there is
no agreed definition of heuristics and metaheuristics
in the literature. Some researchers use `heuristics' and `metaheuristics' interchangeably.
However, the recent trend tends to name all stochastic algorithms
with randomisation and global exploration as metaheuristic.
In this review, we will also follow this convention.

Randomisation provides a good way to move away from
local search to the search on the global scale. Therefore,
almost all metaheuristic algorithms intend to be suitable
for global optimisation.

Metaheuristics can be an efficient way to produce acceptable
solutions by trial and error to a complex problem in
a reasonably practical time.
The complexity of the problem of interest makes it impossible
to search every possible solution or combination, the aim is
to find good feasible solution in an acceptable timescale. There is
no guarantee that the best solutions can be found, and we even do not
know whether an algorithm will work and why if it does work. The idea is
to have an efficient but practical algorithm that will work most
the time and is able to produce good quality solutions. Among the
found quality solutions, it is expected some of them are nearly optimal,
though there is no guarantee for such optimality.

The main components of any metaheuristic algorithms are: intensification
and diversification, or exploitation and exploration (Blum and Roli 2003).
Diversification means to generate diverse solutions so as to explore the search space
on the global scale, while intensification means to focus on the
search in a local region by exploiting the information that a current
good solution is found in this region. This is in combination with
with the selection of the best solutions.

The selection of the best ensures
that the solutions will converge to the optimality.
On the other hand, the diversification
via randomisation avoids the solutions being trapped at local optima, while increases the diversity
of the solutions. The good combination of these two major components will usually ensure
that the global optimality is achievable.

Metaheuristic algorithms can be classified in many ways. One way is to classify them
as: population-based and trajectory-based. For example, genetic algorithms are
population-based as they use a set of strings, so is the
particle swarm optimisation (PSO) which uses multiple
agents or particles (Kennedy and Eberhart 1995).
On the other hand, simulated annealing uses a single
agent or solution which moves through the design space or search space
in a piecewise style (Kirkpatrick et al. 1983).

\section{Overview of Metaheuristics}

Throughout history, especially at the
early periods of human history, the main approach to problem-solving
has always been heuristic or metaheuristic -- by trial and error.
Many important discoveries were done by `thinking outside the box',
and often by accident;  that is heuristics. Archimedes's Eureka moment
was a heuristic triumph. In fact, our daily learning experience
(at least as a child) is dominantly heuristic.

Despite its ubiquitous nature, metaheuristics as a scientific method
to problem solving is indeed a modern phenomenon, though
it is difficult to pinpoint when the metaheuristic method was first used.
Alan Turing was probably the first to use heuristic algorithms
for code-breaking during the Second World War (Copeland 2004).
Turing called his search method {\it heuristic search},
as it could be expected it worked most of time, but there was no
guarantee to find the correct solution, but it was a tremendous success.

The 1960s and 1970s were the two important decades for the development of
evolutionary algorithms. First, John Holland and his collaborators  at the University of Michigan
developed the genetic algorithms in 1960s and 1970s (Holland 1975).
As early as 1962, Holland studied the adaptive system and was the first to use
crossover and recombination manipulations for modeling such system. His seminal book
summarizing the development of genetic algorithms was published in 1975.
In the same year, De Jong finished his important dissertation showing the
potential and power of genetic algorithms for a wide range of objective functions, either
noisy, multimodal or even discontinuous (De Jong 1975).

Briefly speaking, a genetic algorithm (GA) is a search method based on the abstraction of
Darwinian evolution and natural selection of biological systems and representing them in the
mathematical operators: crossover or recombination,
mutation, fitness, and selection of the fittest. Ever since, genetic algorithms have
become so successful in solving a wide range of optimisation problems, there have
several thousands of research articles and hundreds of books written.
Some statistics show that a vast majority of Fortune 500 companies are now using them routinely
to solve tough combinatorial optimisation problems such as planning, data-fitting,
and scheduling.

During the about same period in the 1960s, Ingo Rechenberg and Hans-Paul Schwefel
 both then at the Technical University of Berlin
developed a search technique for solving optimisation
problem in aerospace engineering, called evolutionary strategy, in 1963.
Later, Peter Bienert joined them and began to construct an automatic experimenter
using simple rules of mutation and selection.  There was
no crossover in this technique, only mutation was used to produce an offspring
and an improved solution was kept at each generation. This was essentially
a simple trajectory-style hill-climbing algorithm with randomisation. As early as
1960, Lawrence  J. Fogel intended to use simulated evolution as a learning process
as a tool to study artificial intelligence.  Then,
in 1966,  L. J.  Fogel, together A. J. Owen and M. J. Walsh,
developed the evolutionary programming technique by representing solutions
as finite-state machines and randomly mutating one of these machines (Fogel et al. 1966).
The above innovative ideas and methods have evolved into a much wider
discipline, called evolutionary algorithms and/or evolutionary computation.

A loosely related topic is the development of the support vector machine which is
a classification technique dated back to the earlier work by V. Vapnik in 1963 on linear classifiers,
and the nonlinear classification with kernel techniques were
developed by V. Vapnik and his collaborators in the 1990s
(Vapnik 1995, Vapnik et al. 1997).

A brief history of heuristics till 1960s
in the context of combinatorial optimisation
was reviewed by Schrijver (2005). One of the first books on heuristics
was published by Judea (1984). A very simplified, informal history can be
found from the wikipedia article on metaheuristic (Wikipedia 2010).

The two decades of 1980s and 1990s were the most exciting time for
metaheuristic algorithms.
The next big step is the development of simulated annealing (SA) in 1983, an optimisation
technique, pioneered by S. Kirkpatrick, C. D. Gellat and M. P. Vecchi,
inspired by the annealing process of metals. It is a trajectory-based search
algorithm starting with an initial guess solution at a high temperature,
and gradually cooling down the system (Kirkpatrick et al. 1983).
A move or new solution is accepted
if it is better; otherwise, it is accepted with a probability, which makes
it possible for the system to escape any local optima. It is then expected that
if the system is cooled down slowly enough, the global optimal solution can be reached.

Most metaheuristic do not explicitly use memory, except the selection
of the best solutions. The actual first usage of memory in modern metaheuristics is probably due to Fred Glover's
Tabu search in 1986, though his seminal book on Tabu search was published
later in 1997 (Glover and Luguna 1997).

Interestingly, an even-higher level algorithm, namely, the memetic algorithm, proposed by P. Moscato
in 1989, is a multi-generation, co-evolution and self-generation algorithm (Moscato 1989), and it
can be considered as a hyper-heuristic algorithm, rather than metaheuristic.

At the beginning of the 1990s,
Marco Dorigo finished his PhD thesis on optimisation and natural
algorithms (Dorigo 1992), in which he described his innovative work on
ant colony optimisation (ACO).
This search technique was inspired by the swarm intelligence of social
ants using pheromone as a chemical messenger.
Then, in 1992, John R. Koza of Stanford University published a treatise
on genetic programming which laid  the foundation of a whole new area of machine learning,
revolutionizing computer programming (Koza 1992).
As early as in 1988, Koza applied his
first patent on genetic programming. The basic idea is to use the genetic principle to breed
computer programs so as to gradually produce the best programs for a given type of problem.

Another interesting algorithm is the so-called
artificial immune system,  inspired by the characteristics
of the immune system of mammals to use memory and
learning as a novel approach to problem solving.
The idea was proposed by Farmer et al. in 1986 (Farmer et al. 1986),
with important work on immune networks by Bersini and Varela in 1990 (Bersini and Varela 1990).
It is an adaptive system with high potential. There are many variants
developed over the last two decades, including the clonal selection algorithm,
negative selection algorithm, immune networks and others.

Slightly later in 1995, another significant progress is the development of
the particle swarm optimisation (PSO) by American social psychologist James Kennedy,
and engineer Russell C. Eberhart (Kennedy and Eberhart 1995). Loosely speaking, PSO is an optimisation
algorithm inspired by swarm  intelligence of fish and birds and by
even human behavior. The multiple agents, called particles, swarm around the
search space starting from some initial random guess. The swarm communicates the current best
and shares the global best so as to focus on the quality solutions.

Since the development of PSO, there have been more than about 20 different variants
of particle swarm optimisation techniques, and have been applied
to almost all areas of tough optimisation problem. A variant typically extends
PSO by varying some aspects of the particle swarm or introducing new functionality,
for example, new variants were developed by introducing velocity variations
(Cui and Zeng 2005, Cui and Cai 2009). There is some strong
evidence that PSO is better than traditional search algorithms and even better
than genetic algorithms for many types of problems, though this is far from conclusive.

Around the same time in 1996 and later 1997, R. Storn and K. Price
developed their vector-based evolutionary algorithm, called differential evolution (DE),
and this algorithm proves more efficient than genetic algorithms in many applications (Storn and Price 1997).

Another interesting method is the cross-entropy method developed by Rubinstein in 1997.
It is a generalized Monte Carlo method, based on the rare event simulations.
This algorithm consists of two phases: generation of random samples and
update of the parameters. Here the aim is to minimise the cross entropy (Rubinstein 1997).

In celebrating this series of success, a significant event is the publication
of the `no free lunch theorems' for optimisation in 1997 by
D. H. Wolpert and W. G. Macready, which sent out a shock way to the optimisation
community (Wolpert and Macready 1997).
Researchers have been always trying to find better algorithms, or even
universally robust algorithms, for optimisation, especially for tough NP-hard
optimisation problems. However, these theorems state that if algorithm A performs
better than algorithm B for some optimisation functions, then B will outperform
A for other functions. That is to say, if averaged over all possible
function space, both algorithms A and B will perform on average equally well.
Alternatively, there is no universally better algorithms exist.
That is disappointing, right? Then, people realized that we do not need
the average over all possible functions  for a given optimisation problem.
What we want is to find the best solutions, which has nothing to do with average over
all possible function space. In addition, we can accept the fact that
there is no universal or magical tool, but we do know from our experience
that some algorithms indeed outperform others for given types of optimisation
problems. So the research now focuses on finding the best and most efficient
algorithm(s) for a given problem. The objective is to design better algorithms
for most types of problems, not for all the problems. Therefore, the search is still on.

At the turn of the 21st century, things became even more exciting. First, Zong Woo Geem
et al. in 2001 developed the harmony search (HS) algorithm (Geem et al. 2001),
which has been widely applied
in solving various optimisation problems such as water distribution, transport modelling
and scheduling. This is almost immediately followed by the development of another
algorithm, called the bacterial foraging optimisation by
K. M. Passino in around 2002, inspired by the social foraging behaviour of bacteria
such as {\it Escherichia coli} (Passino 2002).

In 2004, S. Nakrani and C. Tovey proposed the honey bee algorithm
and its application for optimizing Internet hosting centers (Nakrani and Trovey 2004),
which followed by the development of a novel bee algorithm by D. T. Pham et al. in 2005 (Pham et al. 2005)
and the artificial bee colony (ABC) by D. Karaboga in 2005 (Karaboga 2005).
In the same year, a glowworm algorithm was developed for the detection of
multiple source locations applied to collective robots
(Krishnanand and Ghose 2005). Later in 2006, a honey-bees mating optimisation
algorithm was developed by Haddad et al. (2006).

In 2008, Mucherino and Seref proposed monkey search, based on the foraging
behaviour of monkeys (Mucherino and Seref 2008).
In the same year, the author of this article
developed the firefly algorithm (FA) (Yang 2008, Yang 2009, Yang 2010a).
Quite a few research articles on the firefly algorithm then followed,
and this algorithm has attracted a wide range of interests.
In 2009, Xin-She Yang at Cambridge University, UK, and Suash Deb at Raman College of Engineering, India,
introduced an efficient cuckoo search (CS) algorithm (Yang and Deb 2009, Yang and Deb 2010a),
and it has been demonstrated that
CS is far more effective than most existing metaheuristic algorithms including particle
swarm optimisation.  In 2010, I developed a bat-inspired algorithm
for continuous optimisation (Yang 2010b), and its efficiency is quite promising. At about the same,
Yang and Deb also proposed the eagle strategy as a two-stage optimisation strategy (Yang and Deb 2010b).

\section{Characteristics of Metaheuristics}

The efficiency of metaheuristic algorithms can be attributed to the fact
that they imitate the best features in nature, especially the selection of
the fittest in biological systems which have evolved
by natural selection over millions of years.

Two important characteristics of metaheuristics
are: intensification and diversification (Blum and Roli 2003).
Intensification intends to search locally and more intensively,
while diversification makes sure the algorithm explores the search space globally
(hopefully also efficiently).

Furthermore, intensification is also called exploitation,
as it typically searches around the current best solutions
and selects the best candidates or solutions. Similarly,
diversification is also called exploration, as it strives to explore
the search space more efficiently, often by large-scale randomisation.

The fine balance between these two components is very
important to the overall efficiency and performance of an algorithm.
Too little exploration and too much exploitation could cause the
system to be trapped in local optima, which makes it very difficult
or even impossible to find the global optimum. On the other hand,
if too much exploration but too little exploitation, it may be difficult
for the system to converge and thus slows down the overall search performance.
The proper balance itself is an optimisation problem,
and one of the main tasks of designing new algorithms is
to find a certain balance concerning this optimality and/or tradeoff.

Obviously, simple exploitation and exploration are not enough.
During the search, we have to use a proper mechanism or
criterion to select the best solutions. The  most common criterion
is to use the {\it Survival of the Fittest}, that is to keep updating
the current best found so far. In addition, certain elitism is often
used, and this is to ensure the best or fittest solutions are not
lost, and should be passed onto the next generations.

\section{Importance of Randomisation}

There are many ways of carrying out intensification and diversification.
In fact, each algorithm and its variants use different ways of achieving the
balance of between exploration and exploitation.

By analyzing all the metaheuristic algorithms, we can categorically say that
the way to achieve exploration or diversification is mainly by
certain randomisation in combination with a deterministic procedure.
This ensures that the newly generated solutions distribute as diversely as
possible in the feasible search space. One of simplest and yet most commonly
used randomisation techniques is to use
\be \x_{\rm new} ={\bf L}+({\bf U}-{\bf L}) *\epsilon_u, \ee
where ${\bf L}$ and ${\bf U}$ are the lower bound and upper bound, respectively.
$\epsilon_u$ is a uniformly distributed random variable in [0,1].
This is often used in many algorithms such as harmony search,
particle swarm optimisation and firefly algorithm.
It is worth pointing that the use of a uniform distribution is not the only
way to achieve randomisation. In fact, random walks such as
L\'evy flights on a global scale are more efficient.

A more elaborate way to obtain diversification is to use mutation and crossover.
Mutation makes sure new solutions are as far/different as possible,
from their parents or existing solutions; while crossover limits the
degree of over diversification,
as new solutions are generated by swapping parts of the existing solutions.

The main way to achieve the exploitation is to generate new solutions around
a promising or better solution  locally and more intensively. This can be easily achieved
by a local random walk
\be \x_{\rm new}=\x_{\rm old}+ s \; {\bf w}, \label{div-equ-100} \ee
where ${\bf w}$ is typically drawn from a Gaussian distribution with zero mean.
Here $s$ is the step size of the random walk. In general, the step size
should be small enough so that only local neighbourhood is visited.
If $s$ is too large, the region visited can be  too far away from the region of interest, which
will increase the diversification significantly but reduce the intensification
greatly. Therefore, a proper step size should be much smaller than
(and be linked with) the scale of  the problem. For example,
the pitch adjustment in harmony search  and the move in simulated annealing
are a random walk.

If we want to increase the efficiency of this random walk (and thus increase the
efficiency of exploration as well), we can use other forms of random walks such as
L\'evy flights where $s$ is drawn from a L\'evy distribution with large step sizes.
In fact, any distribution with a long tail will help to increase the step size
and distance of such random walks.

Even with the standard random walk, we can use a more selective or controlled
walk around the current best $\x_{\rm best}$, rather than any good solution.
This is equivalent to replacing the above equation by
\be \x_{\rm new} =\x_{\rm best} + s \; {\bf w}. \ee

Some intensification technique is not easy to decode, but may be
equally effective. The crossover operator in evolutionary algorithms is a good example, as it uses
the solutions/strings from parents to form offsprings or new solutions.

In many algorithms, there is no clear distinction or explicit differentiation between intensification
and diversification. These two steps are often intertwined and
interactive, which may, in some cases, become an advantage.
Good examples of such interaction is the genetic algorithms (Holland 1975),
harmony search (Geem et al. 2001)
and bat algorithm (Yang 2010b). Readers can analyze any chosen algorithm to see how
these components are implemented.

In addition, the selection of the best solutions is a crucial component for the
success of an algorithm. Simple, blind exploration and exploitation may not be
effective without the proper selection of the solutions of good quality.
Simply choosing the best may be effective for optimisation problems with
a unique global optimum. Elitism and keeping the best solutions are efficient
for multimodal and multi-objective problems. Elitism in genetic algorithms
and selection of harmonics are good examples of the selection of the fittest.

In contrast with  the selection of the best solutions, an efficient metaheuristic
algorithm should have a way to discard the worse solutions so as to increase
the overall quality of the populations during evolution. This is often achieved by
some form of randomisation and probabilistic selection criteria. For example,
mutation in genetic algorithms acts a way to do this. Similarly, in the cuckoo search,
the castaway of a nest/solution is another good example (Yang and Deb 2009).

Another important issue is the randomness reduction. Randomisation is mainly
used to explore the search space diversely on the global scale, and also, to some extent,
the exploitation  on a local scale. As better solutions are found,
and as the system converges, the degree of randomness should be reduced;
otherwise, it will slow down the convergence. For example, in particle swarm optimisation,
randomness is automatically reduced as the particles swarm
together (Kennedy and Eberhart 1995),
this is because the distance between each particle and the current
global best is becoming smaller and smaller.

In other algorithms, randomness is not reduced and but controlled and selected.
For example, the mutation rate is usually small so as to limit the randomness,
while in simulated annealing, the randomness during iterations may remain the same, but the
solutions or moves are selected and acceptance probability becomes smaller.

Finally,  from the implementation point of view, the actual implementation does
vary, even though the pseudo code should give a good guide and should not
in principle lead to ambiguity. However, in practice, the actual way of
implementing the algorithm does affect the performance to some degree. Therefore,
validation and testing of any algorithm implementation are important (Talbi 2009).

\section{The Generalized Evolutionary Walk Algorithm (GEWA) }

From the above discussion of all the major components and their characteristics,
we realized that a good combination of local search and global search with a proper
selection mechanism should produce a good metaheuristic algorithm, whatever
the name it may be called.

In principle, the global search should be carried out more frequently at the initial
stage of the search or iterations. Once a number of good quality solutions are found,
exploration should be sparse on the global scale, but frequent enough
so as to escape any local trap if necessary. On the other hand,
the local search should be carried out as efficient as possible, so a good local search
method should be used. The proper balance of these two is paramount.

Using these basic components, we can now design a new, generic, metaheuristic
algorithm for optimisation, we can call it the Generalized Evolutional Walk Algorithm (GEWA).
Evolutionary walk is a random walk, but with a biased selection
towards optimality. This is a generalized framework for global optimisation.

There are three major components in our proposed algorithm: 1) global exploration by
randomisation, 2) intensive local search by random walk, and 3) the selection
of the best with some elitism.
The pseudo code of GEWA is shown in Fig. 1.
The random walk should be carried out around the current global best ${\bf g}_*$
so as to exploit the system information such as the current best more effectively. We have
\be \x_{t+1}= {\bf g}_* +{\bf w}, \ee
and
\be {\bf w}=\varepsilon {\bf d}, \ee
where $\varepsilon$ is drawn from a Gaussian distribution or normal distribution
N$(0, \sigma^2)$, and ${\bf d}$ is the step length vector which should be
related to the actual scales of independent variables. For simplicity, we can
take $\sigma=1$.

\begin{figure}
\begin{center}
\begin{minipage}[c]{0.9 \textwidth}
\hrule \vspace{5pt}
\indent  Initialize a population of $n$ walkers $\x_i \; (i=1,2,...,n)$; \\
\indent Evaluate fitness $F_i$ of $n$ walkers; \\
\indent Find the current best ${\bf g}_*$; \\
\indent {\bf while} ($t<$MaxGeneration) or (stop criterion); \\
\indent \quad Discard worst solutions and replace by (\ref{gewa-eq-550}) or (\ref{gewa-eq-555}); \\
\indent \quad  {\bf if } (rand $<\alpha$), \\
\indent \qquad Local search: random walk around the best
\be \quad \x_{t+1}={\bf g}_*+\varepsilon {\bf d} \label{gewa-eq-550} \ee
\indent \quad  {\bf else} \\
\indent \qquad Global search: randomisation \be \quad \x_{t+1}={\bf L}+({\bf U}-{\bf L}) \; \epsilon_u \label{gewa-eq-555} \ee
\indent \quad {\bf end} \\
\indent \quad Evaluate new solutions \& find the current best ${\bf g}_*^t$; \\
\indent \quad $t=t+1$; \\
\indent  {\bf end while} \\
\indent  Postprocess results and visualization;
\hrule \vspace{5pt}
\caption{Generalized Evolutionary Walk Algorithm. \label{gewa-fig-100} }
\end{minipage}
\end{center}
\end{figure}

The randomisation step can be achieved by
\be \x_{t+1}={\bf L} + ({\bf  U}-{\bf L}) \epsilon_u, \ee
where $\epsilon_u$ is drawn from a uniform distribution Unif[0,1].
${\bf U}$ and ${\bf L}$ are the upper and lower bound vectors, respectively.

There are only two algorithm-dependent parameters in GEWA: the population size $n$
and the randomization control parameter $\alpha$.
Typically, $\alpha \approx 0.25 \sim 0.7$, and we often use $\alpha=0.5$ in
most applications. Interested readers can try to do some parametric studies.
Furthermore, the number ($n$) of random walkers is also important. Too few walkers
are not efficient, while too many may lead to slow convergence.
In general, the choice of $n$ should follow the similar guidelines
as those for all population-based algorithms. Typically, we can use
$n=15$ to $50$ for most applications.

Again two important issues are: 1) the balance of intensification and diversification
controlled by a single parameter $\alpha$, and 2) the choice of the
step size of the random walk. Parameter $\alpha$ is typically
in the range of $0.25$ to $0.7$.
The choice of the right step size is also important, and
the ratio of the step size to its length scale
can be typically $0.001$ to $0.01$ for most applications.
Another important issue is the selection of the best and/or elitism.
As we intend to discard the worst solution and replace it by generating
new solution. This may implicitly weed out the least-fit solutions, while
the solution with the highest fitness remains in the population.
The selection of the best and elitism can be guaranteed implicitly
in the evolutionary walkers.

Preliminary studies have shown that GEWA is as efficient as many other
metaheuristic algorithms, and its working mechanism is simple.
A detailed analysis of such comparison and analysis
will be reported elsewhere in another paper.

\section{Open Problems}

Despite the success of modern metaheuristic algorithms,
there are many important
questions which remain unanswered. We know how these heuristic algorithms work, and we also partly
understand why these algorithms work. However,  it is difficult
to analyze mathematically why these algorithms are so successful. In fact, these
are unresolved open problems.

Apart from the mathematical analysis on simulated annealing and particle swarm optimisation,
 convergence of all other algorithms has not been proved mathematically, at least
up to now. Any mathematical analysis will thus provide important insight into these algorithms. It will also be
valuable for providing new directions for further important modifications on these algorithms
or even pointing out innovative ways of developing new algorithms.

In addition, it is still only partly understood why different components of heuristics and metaheuristics
interact in a coherent and balanced way so that they produce efficient algorithms which converge under the given
conditions. For example, why does a balanced combination of randomisation and a deterministic component
lead to a much more efficient algorithm (than a purely deterministic and/or a purely random algorithm)?
How to measure or test if a balance is reached? How to prove that the use of memory can significantly increase the
search efficiency of an algorithm? Under what conditions?

From the well-known No-Free-Lunch theorems (Wolpert and Macready 1997),
 we know that they have been proved for single objective optimisation,
but they remain unproved for multiobjective optimisation. If they are proved to be true (or not)
for multiobjective optimisation, what are the implications for algorithm development?

Furthermore, there is no agreed measure for comparing performance of different algorithms, thought
the absolute objective value and the number of functional evaluations are two widely used
measures. An attempt to provide a general framework was given by Shilane et al. (2008).
However, a formal theoretical analysis is yet to be developed.
If you are looking for some research topics, either for yourself or for your research students,
these could form important topics for further research.

Another interesting question that people often ask me is ``to be inspired or not to be inspired"?
Nature has evolved over billions of years, she has found almost perfect solutions
to every problem she has met. Almost all the not-so-good solutions have been discarded
via natural selection. The  optimal solutions seem (often after a huge number of
generations) to appear at the evolutionarilly stable equilibrium, even though we may
not understand how the perfect solutions are reached.
When we try to solve human problems, why not try to be inspired by
the nature's success? The simple answer to the
question `To be inspired or not to be inspired?' is
`why not?'. If we do not have good solutions at hand, it is always
a good idea to learn from nature.

Nature provides almost unlimited ways for problem-solving. If we can observe carefully,
we are surely inspired to develop more powerful and efficient new generation algorithms.
Intelligence is a product of biological evolution in nature.
Ultimately some intelligent algorithms (or systems) may appear in the future,
so that they can evolve and optimally adapt to solve
NP-hard optimisation problems efficiently and intelligently.

In addition, a current trend is to use simplified metaheuristic algorithms to
deal with complex optimisation problems. Possibly, there is a need to
develop more complex metaheuristic algorithms which can truly mimic
the exact working mechanism of some natural or biological systems,
leading to more powerful next generation, self-regulating, self-evolving,
and truly intelligent metaheuristics.

\end{document}